\documentstyle{article}
\newtheorem{theorem}{Theorem}[section]

\newcommand{\si}{\sigma}

\newcommand{\mbf}{\mathbf}
\newcommand{\dsty}{\displaystyle}
\newcommand{\ra}{\rightarrow}

\newcommand{\al}{\alpha}
\newcommand{\be}{\beta}
\newcommand{\de}{\delta}

\newcommand{\ga}{\gamma}
\newcommand{\la}{\lambda}

\begin{document}

\title{ Random Energy Model with Compact Distributions}
\author{Nabin Kumar Jana\footnote{Stat-Math Unit, Indian Statistical Institute, 203 B. T. Road, Kolkata,
India.
e-mail: nabin\_r@isical.ac.in}\\
 {\small Indian Statistical Institute, Kolkata, India}}

\date{February 8, 2005}
\pagestyle{myheadings} \markboth{Jana N. K.}{Random Energy Model}

\maketitle
\begin{abstract}
In this paper we study the Random energy model - so called toy model
of the spin glass theory - where the underlying distributions are
compactly supported. We prove a general theorem on the asymptotics
of free energy and obtain formulae in several interesting cases -
like uniform distribution, truncated double
exponential.\\

\end{abstract}

\section{Introduction}
Random Energy Model (REM) is a simple model in the theory of spin
glasses, originally proposed by Derrida$^{\cite{D}}$. For each
configuration $\si \in\{-1,1\}^N$ of an $N$ particle system, the
Hamiltonian $H_N(\si)$ is random and they are i.i.d. over $\si$.
Apart from the study of distribution of the Gibbs' distribution, one
of the problems in REM is the study of asymptotics of ${\mbf{E}}\log
Z_N(\be)$, where $Z_N(\be)$ is the partition function defined as
$\sum_\si e^{-\be H_N(\si)}$ for $\be>0$. This is well studied in
the literature$^{\cite{B,CEGG,DW,T1}}$ when $H_N$ are Gaussian.
Guerra's$^{\cite{GT}}$ convexity arguments prove the existence of
$\dsty\lim_{N\ra \infty}\frac{1}{N}{\mbf{E}}\log Z_N(\be)$ even when
$H_N(\si)$ are non Gaussian provide a moment condition holds.
Guerra's treatment was for the SK-Model. See Contucci {\em et
al}$^{\cite{CEGG}}$ where this is done for REM, though we have
trouble following their convexity argument for the Generalized
Random Energy Model. The methods - powerful as they are - will not
help in the evaluation of the limits. In this paper we prove a
general theorem when the support of distribution of $H_N(\si)$ is
compact. This exercise is undertaken for several reasons. The
powerful exponential inequalities are not in general available in
the non-Gaussian case. More so, when the distribution is flat, like
uniform. The rate at which length of the supports grow should have
some influence on the asymptotics. The main theorem specifies
certain conditions when there is no phase transition. Our arguments
for lower bound are similar in spirit to
Talagrand$^{\cite{T1}}$(P.11-12).

\section{Main Result}
For each $N$, let $H_N(\si)$ for $\si \in\{-1,1\}^N$ be $2^N$ i.i.d.
symmetric random variables having density $\phi_N$ with compact
support. Let $[-T_N,T_N]$ be the support of $\phi_N$. As expected,
large value of $H_N(\si)$ are the most relevant. Accordingly, for
$0< s <1$, define $a_N(s)= {\mbf{P}}\{H_N(\si) \geq sT_N\}$.

\begin{theorem}
Suppose for any $s$, $0<s<1$, there exists $m\geq 0$ (possibly
depending on $s$) such that $0< \dsty\lim_{N\ra\infty} N^m a_N(s)
<\infty$ and $\dsty\lim_{N\ra \infty} \frac{N}{T_N} \ra \al$, $0
\leq \al \leq \infty$. Then for $\be>0$
$$\dsty\lim_{N\ra \infty} \frac{1}{T_N}{\mbf{E}}\log Z_N(\be) = \al
\log 2 + \be$$  (when $\al = \infty$, the right side is interpreted
as $\infty$).
\end{theorem}

\noindent {\bf Remark 2.1} Of course, the result can equivalently be
stated as
$$\dsty\lim_{N\ra \infty} \frac{1}{N}{\mbf{E}} \log Z_N(\be) = \log
2 + \frac{\be}{\al}.$$

\noindent {\it Proof:} First consider the case $\al < \infty$.

Since $\log$ is concave, by Jensen's inequality
\begin{equation}
{\mbf{E}}\log Z_N(\be) \leq \log {\mbf{E}}Z_N(\be).
\end{equation}

As $H_N$'s are symmetric and i.i.d, $${\mbf{E}}Z_N(\be) = 2^N
{\mbf{E}}e^{\be H_N} < 2^N e^{\be T_N}.$$

Hence, by assumption and (1),
\begin{equation}
\dsty\limsup_{N \ra \infty}\frac{1}{T_N}{\mbf{E}}\log Z_N(\be) \leq
\al \log 2 + \be.
\end{equation}

We now show,
\begin{equation}
\dsty\liminf_{N \ra \infty}\frac{1}{T_N}{\mbf{E}}\log Z_N(\be) \geq
\al \log 2 + \be.
\end{equation}
Fix $0<s<1$ and let $X_N= \#\{\si: -H_N(\si) \geq s T_N\}$. Then
${\mbf{E}}X_N = 2^N a_N(s)$ and ${\mbf{E}}X_N^2 = 2^N(2^N-1)a_N^2(s)
+ 2^N a_N(s)$, so that
\begin{equation}
{\mbf{E}}(X_N-{\mbf{E}}X_N)^2 \leq 2^N a_N(s).
\end{equation}

If $A_N = \{X_N \leq 2^{N-1} a_N(s)\}$ then $A_N \subset
\{(X_N-{\mbf{E}}X_N)^2 \geq 2^{2N-2} a_N^2(s)\}$, so that by Markov
inequality and (4), $${\mbf{P}}(A_N) \leq
\frac{{\mbf{E}}(X_N-{\mbf{E}}X_N)^2}{2^{2N-2}a_N^2(s)} \leq
\frac{4}{2^Na_N(s)},$$

i.e., ${\mbf{P}}(A_N^c) \geq 1- \frac{4}{2^Na_N(s)}$. But on
$A_N^c$,
$$Z_N(\be) \geq X_N e^{\be s T_N} > 2^{N-1} a_N(s) e^{\be s
T_N},$$ and hence
\begin{equation}
{\mbf{E}}\log Z_N(\be) {\mbf{1}}_{A_N^c} > [(N-1)\log2 + \log a_N(s)
+ \be s T_N](1-\frac{4}{2^N a_N(s)}).
\end{equation}

Since always, $Z_N(\be) \geq 2^N e^{-\be T_N}$ and ${\mbf{P}}(X_N=0)
= (1-a_N(s))^{2^N}$, we have
\begin{equation}
{\mbf{E}}\log Z_N(\be){\mbf{1}}_{\{X_N=0\}} \geq (N\log2 -\be
T_N)(1-a_N(s))^{2^N}.
\end{equation}

On $\{1 \leq X_N \leq 2^{N-1}a_N(s)\}$, $\log Z_N(\be)
>\dsty\max_\si \{-H_N(\si)\} \geq \be s T_N >0$ and hence
\begin{equation}
{\mbf{E}}\log Z_N(\be){\mbf{1}}_{\{1 \leq X_N \leq 2^{N-1}a_N(s)\}}
> 0.
\end{equation}

Since $A_N= \{X_N=0\}\cup\{1\leq X_N \leq 2^{N-1} a_n(s)\}$ using
(5), (6) and (7) we have
$$\begin{array}{lll}
\frac{1}{T_N}{\mbf{E}}\log Z_N(\be) &\geq & [\frac{N-1}{T_N}\log2+
\frac{\log a_N(s)}{T_N} +\be s](1-\frac{4}{2^Na_N(s)})\\
&& \hspace{18ex} + [\frac{N}{T_N}\log2 -\be](1-a_N(s))^{2^N}.
\end{array}$$

By assumptions,
$$\frac{\log a_N(s)}{T_N} = \frac{1}{T_N} \log(N^m
a_N(s)) - m\frac{\log N}{T_N} \ra 0$$

and $2^Na_N(s) \ra \infty$ so that $(1-a_N(s))^{2^N} \ra 0$ as $N\ra
\infty$. Thus, under the assumptions,
$$\dsty\liminf_{N\ra \infty}\frac{1}{T_N}{\mbf{E}}\log Z_N(\be)
\geq \al \log2 + \be s.$$

This being true for any $0<s<1$, (3) follows.

If $\al=\infty$ the above argument shows that
$$\dsty\liminf_{N\ra \infty}\frac{1}{T_N}{\mbf{E}}\log Z_N(\be)= \infty$$
to complete the proof.

\section{Examples}
\noindent {\bf Example 3.1}({\em Truncated Double Exponential})

Let $$\phi_N(x) =\frac{\la_N}{2(1-e^{-\la_N T_N})}e^{-\la_N
|x|}{\mbf{1}}_{[-T_N,T_N]}.$$ Then for each $s$, $0<s<1$, we have
$\dsty\lim_{N\ra \infty}a_N(s)>0$ if $\la_N T_N \ra \de$($0\leq \de
< \infty$) as $N\ra \infty$. Thus, for example, if $T_N =N$ and
$\la_N= \frac{\de}{N} (0\leq \de <\infty)$ then $$\dsty\lim_{N\ra
\infty} \frac{1}{N}{\mbf{E}}\log Z_N(\be) = \log2 + \be,$$ where as,
if $T_N =N^2$ and $\la_N=\frac{\de}{N^2} (0\leq \de <\infty)$ then
$$\dsty\lim_{N\ra \infty} \frac{1}{N^2}{\mbf{E}}\log Z_N(\be) = \be$$

\noindent {\bf Example 3.2}({\em Uniform Distribution})

Let $$\phi_N(x) =\frac{1}{2T_N}{\mbf{1}}_{[-T_N,T_N]}.$$ Then
$a_N(s) = \frac{1-s}{2} >0$ for all $0<s<1$. The theorem now
implies,

a) if $T_N = N^\ga$ with $0<\ga<1$ or $T_N=\log N$ then
$\dsty\lim_{N\ra \infty} \frac{1}{T_N}{\mbf{E}}\log Z_N(\be) =
\infty$,

b) if $T_N = N$ then $\dsty\lim_{N\ra \infty}
\frac{1}{N}{\mbf{E}}\log Z_N(\be) = \log2 + \be$,

c) if $T_N = N^\ga$ with $\ga>1$ or $T_N=2^N$ then $\dsty\lim_{N\ra
\infty} \frac{1}{T_N}{\mbf{E}}\log Z_N(\be) = \be$.

Similar remarks apply for the following examples also.\\

\noindent {\bf Example 3.3}({\em Bernoulli})

Let $$\phi_N(x) = \frac{1}{2}[\de_{-T_N} + \de_{T_N}].$$ Here
$a_N(s)= \frac{1}{2}$ for all $0<s<1$.\\

\noindent {\bf Example 3.4}

Let $$\phi_N(x) =\frac{\la + 1}{2T_N^{\la+1}}(T_N-|x|)^\la,
\hspace{5ex}-T_N\leq x \leq T_N \mbox{ and } \la >0.$$ Then $a_N(s)
= \frac{1}{2}(1-s)^{\la+1}>0$ for all $0<s<1$.\\

\noindent {\bf Example 3.5}

Let $$\phi_N(x) =\frac{1}{2T_N}\cos\frac{x}{T_N},
\hspace{5ex}-\frac{\pi T_N}{2}\leq x \leq \frac{\pi T_N}{2}.$$ Then
$a_N(s) = \frac{1}{2}(1-\sin\frac{s\pi}{2})>0$ for all $0<s<1$.\\

\noindent {\bf Example 3.6}

Let $$\phi_N(x)
=\frac{N}{2(e^{T_N}-1)}e^{N|x|}{\mbf{1}}_{[-T_N,T_N]}.$$ Then
$a_N(s)\ra \frac{1}{2} >0$ as $N\ra \infty$ for all $0<s<1$.\\

\noindent {\bf Remarks 3.1} We are not clear if the convergence in
the theorem holds almost everywhere instead of in expectation.\\

\noindent{\bf Remarks 3.2} In some cases we have been able to
establish limiting law for Gibbs' distributions, but a general
result has not yet emerged.


\section*{Acknowledgments}

The author would like to thank B. V. Rao for valuable discussions.



\end{document}